\documentclass[11pt]{amsart}
\input{diagrams}

\usepackage{amsmath}
\usepackage{amsthm}
\usepackage{amssymb}
\usepackage{amsfonts}
\usepackage{amsxtra}
\usepackage{amscd}
\usepackage{epsfig}
\usepackage{verbatim}
\usepackage{latexsym,amstext,epsfig}

\usepackage[all, knot]{xy}
\xyoption{arc}

\newsymbol\pp 1275

\newcommand{\Hom}{\operatorname{Hom}}

\newcommand{\Rep}{\operatorname{Rep}}
\newcommand{\Proj}{\operatorname{Proj}}

\newcommand{\SI}{\operatorname{SI}}
\newcommand{\SL}{\operatorname{SL}}
\newcommand{\GL}{\operatorname{GL}}
\newcommand{\PGL}{\operatorname{PGL}}
\newcommand{\ZZ}{\mathbb Z}
\newcommand{\CC}{\mathbb C}
\newcommand{\RR}{\mathbb R}

\newcommand{\QQ}{\mathbb Q}

\newcommand{\Id}{\operatorname{Id}}

\newcommand{\rel}{\operatorname{relint}}

\newtheorem{theorem}{Theorem}[section]
\newtheorem{proposition}[theorem]{Proposition}
\newtheorem{corollary}[theorem]{Corollary}
\newtheorem{lemma}[theorem]{Lemma}

\theoremstyle{definition}
\newtheorem{definition}[theorem]{Definition}
\newtheorem{remark}[theorem]{Remark}

\title[]{Notes on GIT-fans for quivers}

\author{Calin Chindris}
\address{University of Minnesota, School of Mathematics, Minneapolis, MN, USA}
\email[Calin Chindris]{chindris@math.umn.edu}

\date{25 August, 2007; Revised: \today}

\begin{document}
\bibliographystyle{plain}

\maketitle

\section{Introduction}

Given a quiver $Q$ and a dimension vector $\beta$, King's various
moduli spaces of $\beta$-dimensional semi-stable representations
depend on the choice of a weight. A central question of Geometric
Invariant Theory is to understand how these quotient varieties
vary as the weights vary. This problem of variation of quotients
for quivers was considered by Hille and de la Pena in \cite{HiP}
where they introduced the wall system associated to $(Q,\beta)$.

Our goal in this short note is to go through the construction of
GIT-fans for quivers. Let $C(Q,\beta)$ be the cone of effective
weights associated to $(Q,\beta)$. Then, the cones defining the
GIT-fan have the property that their relative interiors are
precisely the GIT-classes.

\begin{theorem}[The GIT-fan] \label{thm-GIT-fan} Let $Q$ be a quiver without
oriented cycles and $\beta$ a dimension vector. Define the GIT-
cone associated to an effective weight $\sigma \in C(Q,\beta)$ by
$$
C(\sigma)=\{\sigma' \in C(Q,\beta) \mid
\Rep(Q,\beta)^{ss}_{\sigma} \subseteq \Rep(Q,\beta)^{ss}_{\sigma'}
\}.
$$
Then:
\begin{enumerate}
\renewcommand{\theenumi}{\arabic{enumi}}

\item for every effective weight $\sigma,$ $C(\sigma)$ is a
rational convex polyhedral cone whose relative interior is the
GIT-class of $\sigma$;

\item the cones $C(\sigma)$ form a finite fan covering of
$C(Q,\beta)$.
\end{enumerate}
Consequently, there is a finite number of open subsets of
$\Rep(Q,\beta)$ which can be realized as sets of
$\sigma$-semi-stable representations for some effective weights
$\sigma$.
\end{theorem}

The construction of the GIT-fans for normal projective varieties
is due to Ressayre \cite[Theorem 5.2]{Re}. In \cite[Theorem
3.2]{Ha}, Halic used Ressayre's theorem to prove the existence of
GIT-fans in the affine case. Berchtold and Hausen \cite[Theorem
2.11]{BeHa} gave a simple construction of the GIT-fans for torus
actions on affine varieties. This case was recently extended to
arbitrary reductive groups by Arzhantsev and Hausen in
\cite{ArHa2}.

We prove Theorem \ref{thm-GIT-fan} from scratch following closely
the steps outlined in \cite{Re} (but avoiding Dolgachev and Hu's
finiteness theorem \cite[Theorem 1.3.9]{DH}). Our arguments are
based on King \cite{K} semi-stability criterion for quiver
representations and Schofield's \cite{S1} theory of general
representations of quivers. Let us point out that when $\beta$ is
the thin sincere dimension vector (all coordinates equal to one),
the GIT-fan in the theorem above coincides with the fan
constructed by Hille in \cite[Theorem 4.1]{Hi}.











\section{Recollection from quiver invariant theory}
\label{sec.intro}

Let $Q=(Q_0,Q_1,t,h)$ be a finite quiver without oriented cycles,
where $Q_0$ is the set of vertices, $Q_1$ is the set of arrows and
$t,h:Q_1 \to Q_0$ assign to each arrow $a \in Q_1$ its tail
\emph{ta} and head \emph{ha}, respectively.

For simplicity, we will be working over the field of complex
numbers $\CC$. A representation $V$ of $Q$ over $\CC$ is a family
of finite dimensional $\CC$-vector spaces $\lbrace V(x) \mid x\in
Q_0\rbrace$ together with a family $\{ V(a):V(ta)\rightarrow V(ha)
\mid a \in Q_1 \}$ of $\CC$-linear maps. If $V$ is a
representation of $Q$, we define its dimension vector $\underline
d_V$ by $\underline d_V(x)=\dim_{\CC} V(x)$ for every $x\in Q_0$.
Thus the dimension vectors of representations of $Q$ lie in
$\Gamma=\ZZ^{Q_0}$, the set of all integer-valued functions on
$Q_0$.

Given two representations $V$ and $W$ of $Q$, we define a morphism
$\phi:V \rightarrow W$ to be a collection of linear maps $\lbrace
\phi(x):V(x)\rightarrow W(x)\mid x \in Q_0 \rbrace$ such that for
every arrow $a\in Q_1$, we have $\phi(ha)V(a)=W(a)\phi(ta)$. We
denote by $\Hom_Q(V,W)$ the $\CC$-vector space of all morphisms
from $V$ to $W$. In this way, we obtain the abelian category
$\Rep(Q)$ of all quiver representations of $Q.$ Let $W$ and $V$ be
two representations of $Q$. We say that $V$ is a subrepresentation
of $W$ if $V(x)$ is a subspace of $W(x)$ for all vertices $x \in
Q_0$ and $V(a)$ is the restriction of $W(a)$ to $V(ta)$ for all
arrows $a \in Q_1$.

Let $\beta \in \ZZ_{\geq 0}^{Q_0}$ be a dimension vector of $Q.$
We write $\beta_1\hookrightarrow \beta$ if every $\beta$-
dimensional representation has a subrepresentation of dimension
vector $\beta_1$. If $\sigma \in \mathbb R^{Q_0}$ and $\beta \in
\ZZ^{Q_0}$, we define $\sigma(\beta)$ by
$$
\sigma(\beta)=\sum_{x \in Q_0}\sigma(x)\beta(x).
$$

The \emph{cone of effective weights} associated to $(Q, \beta)$ is
defined by
$$
C(Q,\beta)= \{ \sigma \in \RR^{Q_0}  \mid \sigma(\beta)=0 \text{
and } \sigma(\beta_1) \leq 0  \text{ for all } \beta_1
\hookrightarrow \beta \}.$$

The (saturated, affine) semigroup of lattice points of
$C(Q,\beta)$ is
$$
\Sigma(Q,\beta)=\{ \sigma \in \ZZ^{Q_0}  \mid \sigma(\beta)=0
\text{ and } \sigma(\beta_1) \leq 0  \text{ for all } \beta_1
\hookrightarrow \beta \}.
$$

\subsection{Semi-invariants and moduli spaces for quivers} For every vertex $x$, we denote by
$\varepsilon_x$ the simple dimension vector corresponding to $x$,
i.e., $\varepsilon_x(y)=\delta_{x,y}, \forall y\in Q_0,$ where
$\delta_{x,y}$ is the Kronecker symbol. If $\alpha,\beta$ are two
elements of $\Gamma$, we define the Euler form by
$$
\langle\alpha,\beta \rangle = \sum_{x \in Q_0}
\alpha(x)\beta(x)-\sum_{a \in Q_1} \alpha(ta)\beta(ha).
$$

Let $\beta$ be a dimension vector of $Q$. The representation space
of $\beta-$dimensional representations of $Q$ is defined by
$$\Rep(Q,\beta)=\bigoplus_{a\in Q_1}\Hom(\CC^{\beta(ta)}, \CC^{\beta(ha)}).$$
If $\GL(\beta)=\prod_{x\in Q_0}\GL(\beta(x))$ then $\GL(\beta)$
acts algebraically on $\Rep(Q,\beta)$ by simultaneous conjugation,
i.e., for $g=(g(x))_{x\in Q_0}\in \GL(\beta)$ and $V=\{V(a)\}_{a
\in Q_1} \in \Rep(Q,\beta)$, we define $g \cdot V$ by
$$(g\cdot V)(a)=g(ha)V(a)g(ta)^{-1}\ \text{for each}\ a \in Q_1.$$ Hence, $\Rep(Q,\beta)$ is a rational representation of the
linearly reductive group $\GL(\beta)$ and the $\GL(\beta)-$orbits
in $\Rep(Q,\beta)$ are in one-to-one correspondence with the
isomorphism classes of $\beta-$dimensional representations of $Q$.
As $Q$ is a quiver without oriented cycles, one can show that
there is only one closed $\GL(\beta)-$orbit in $\Rep(Q,\beta)$ and
hence the invariant ring $\text{I}(Q,\beta)= \CC
[\Rep(Q,\beta)]^{\GL(\beta)}$ is exactly the base field $\CC$.

Now, consider the subgroup $\SL(\beta) \subseteq \GL(\beta)$
defined by
$$
\SL(\beta)=\prod_{x \in Q_0}\SL(\beta(x)).
$$

Although there are only constant $\GL(\beta)-$invariant polynomial
functions on $\Rep(Q,\beta)$, the action of $\SL(\beta)$ on
$\Rep(Q,\beta)$ provides us with a highly non-trivial ring of
semi-invariants. Note that any $\sigma \in \ZZ^{Q_0}$ defines a
rational character of $\GL(\beta)$ by
$$\{g(x) \mid x \in Q_0\} \in \GL(\beta) \mapsto \prod_{x \in
Q_0}(\det g(x))^{\sigma(x)}.$$ In this way, we can identify
$\Gamma=\ZZ ^{Q_0}$ with the group $X^\star(\GL(\beta))$ of
rational characters of $\GL(\beta),$ assuming that $\beta$ is a
sincere dimension vector (i.e., $\beta(x)>0$ for all vertices $x
\in Q_0$). We also refer to the rational characters of
$\GL(\beta)$ as (integral) weights.

Let $\SI(Q,\beta)= \CC [\Rep(Q,\beta)]^{\SL(\beta)}$ be the ring
of semi-invariants. As $\SL(\beta)$ is the commutator subgroup of
$\GL(\beta)$ and $\GL(\beta)$ is linearly reductive, we have
$$\SI(Q,\beta)=\bigoplus_{\sigma
\in X^\star(\GL(\beta))}\SI(Q,\beta)_{\sigma},
$$
where $$\SI(Q,\beta)_{\sigma}=\lbrace f \in \CC [\Rep(Q,\beta)]
\mid g \cdot f= \sigma(g)f \text{~for all~}g \in
\GL(\beta)\rbrace$$ is the space of semi-invariants of weight
$\sigma.$

\begin{definition}\cite[Proposition 3.1]{K} Let $\beta$ be a dimension vector and
$\sigma$ an integral weight such that $\sigma(\beta)=0.$ A
$\beta$-dimensional representation $W$ is said to be:
\begin{enumerate}
\renewcommand{\theenumi}{\arabic{enumi}}

\item \emph{$\sigma$-semi-stable} if $\sigma(\underline d_{W'})
\leq 0$ for every subrepresentation $W'$ of $W;$

\item \emph{$\sigma$-stable} if $\sigma(\underline d_{W'})<0$ for
every proper subrepresentation $0 \neq W' \varsubsetneq W.$
\end{enumerate}
\end{definition}

We say that a dimension vector $\beta$ is
\emph{$\sigma$(-semi)-stable} if there exists a
$\beta$-dimensional representation which is
$\sigma$(-semi)-stable.

Let $\beta$ be a $\sigma$-semi-stable dimension vector. The set of
$\sigma$-semi-stable representations in $\Rep(Q,\beta)$ is denoted
by $\Rep(Q,\beta)^{ss}_{\sigma}$ while the set of $\sigma$-stable
representations in $\Rep(Q,\beta)$ is denoted by
$\Rep(Q,\beta)^{s}_{\sigma}.$ The one dimensional torus
$$
T=\{(t\Id_{\beta(x)})_{x \in Q_0} \mid t \in \CC^{*} \} \subseteq
\GL(\beta)
$$
acts trivially on $\Rep(Q,\beta)$ and so there is a well-defined
action of $\PGL(\beta)={\GL(\beta)/T}$ on $\Rep(Q,\beta).$ Using
methods from geometric invariant theory, one can construct the
following GIT-quotient of $\Rep(Q,\beta)$:
$$
\mathcal M (Q,\beta)^{ss}_{\sigma}=\Proj (\oplus_{n \geq 0}\SI(Q,
\beta)_{n\sigma}).
$$
It was proved by King \cite{K} that $\mathcal M
(Q,\beta)^{ss}_{\sigma}$ is a categorical quotient of
$\Rep(Q,\beta)^{ss}_{\sigma}$ by $\PGL(\beta).$ Note that
$\mathcal M(Q,\beta)^{ss}_{\sigma}$ is an irreducible projective
variety, called the moduli space of $\beta$-dimensional
$\sigma$-semi-stable representations (for more details, see
\cite{K}).

\subsection{More on semi-stability and GIT-classes}

If $\sigma \in \RR^{Q_0}$ is a weight, define the set of
$\beta$-dimensional $\sigma$-semi-stable representations by
$$
\Rep(Q,\beta)^{ss}_{\sigma}=\{W \in \Rep(Q,\beta) \mid
\sigma(\underline{d}_W)=0 \text{~and~} \sigma(\underline{d}_{W'})
\leq 0, \forall W' \leq W \}.
$$
We know that when $\sigma$ is an integral weight, this set is a
(possibly empty) open subset of $\Rep(Q,\beta).$ In fact, this is
going to be true (see Corollary \ref{coro-inte-wt}) for arbitrary
weights.

Given a weight $\sigma \in \RR^{Q_0}$, we can construct the full
subcategory $\Rep(Q)^{ss}_{\sigma}$ of $\Rep(Q)$ consisting of all
$\sigma$-semi-stable representations of $Q$. It turns out that
$\Rep(Q)^{ss}_{\sigma}$ is an abelian category whose simple
objects are precisely the $\sigma$-stable representations.
Furthermore, every object of $\Rep(Q)^{ss}_{\sigma}$ has a
Jordan-H{\"o}lder filtration whose factors are $\sigma$-stable
representations.

Following Dolgachev and Hu \cite{DH}, we say that two (effective)
weights $\sigma_1$ and $\sigma_2$ are \emph{GIT-equivalent}, and
write $\sigma_1 \thicksim \sigma_2$, if
$$
\Rep(Q,\beta)^{ss}_{\sigma_1}=\Rep(Q,\beta)^{ss}_{\sigma_2}.
$$
For example, $\sigma \thicksim m\sigma$ for any positive integer
$m.$ The GIT-equivalence class of a weight $\sigma \in C(Q,\beta)$
is denoted by $\langle \sigma \rangle.$

The \emph{GIT-cone} associated to a weight $\sigma \in \RR^{Q_0}$
is defined by:
$$
C(\sigma)=\{ \sigma' \in \RR^{Q_0} \mid
\Rep(Q,\beta)^{ss}_{\sigma} \subseteq \Rep(Q,\beta)^{ss}_{\sigma'}
\}.
$$
Note that two (effective) weights $\sigma_1$ and $\sigma_2$ are
GIT-equivalent if and only if $\sigma_1 \in C(\sigma_2)$ and
$\sigma_2 \in C(\sigma_1).$

Let us record the following useful finiteness result:

\begin{lemma}\label{finiteness-lemma}
The GIT-cones are rational convex polyhedral cones and there only
finitely many of them.
\end{lemma}

\begin{proof} For every
$\sigma \in C(Q,\beta),$ let $D_{\sigma}$ be the set of all
(sub)dimension vectors of the form $\beta'=\underline{d}_{W'},$
where $W'$ is a subrepresentation of a $\sigma$-semi-stable
representation $W \in \Rep(Q,\beta).$ Then:
$$
C(\sigma)=\{\sigma' \in \RR^{Q_0} \mid \sigma'(\beta)=0 \text{ and
} \sigma'(\beta') \leq 0  \text{ for all } \beta' \in D_{\sigma}
\}.
$$
Since there are finitely many subdimension vectors $\beta' \leq
\beta,$ we clearly have only finitely many GIT-cones.
\end{proof}

\begin{corollary} \label{coro-inte-wt}
Any effective weight is GIT-equivalent to some \emph{integral}
effective weight.
\end{corollary}

\begin{proof}
Let $\sigma \in C(Q,\beta).$ As $C(\sigma)$ is a rational convex
polyhedral cone, we can find a sequence $\{\sigma_n\}_{n\geq 1}
\subseteq C(\sigma)\bigcap \QQ^{Q_0}$ such that $\sigma_n \to
\sigma.$ Since there are finitely many GIT-cones, we can assume
that $C(\sigma_n)=C(\sigma_1),$ $\forall n \geq 1.$ This implies
that $\sigma \in C(\sigma_1)$ as $C(\sigma_1)$ is a closed subset
of $\RR^{Q_0}.$ Hence, $\sigma \thicksim \sigma_1$ and we are
done.
\end{proof}

\begin{definition} Let $\sigma \in C(Q,\beta)$ and $W \in \Rep(Q,\beta)^{ss}_{\sigma}.$
We say that $W$ is \emph{$\sigma$-polystable} if $\GL(\beta)W$ is
closed in $\Rep(Q,\beta)^{ss}_{\sigma}.$
\end{definition}

\begin{remark} Note that the $\sigma$-polystable points
are said to be pivotal for $\sigma$ in the terminology of
\cite{Re}. Our terminology is motivated by the following result of
King:
\end{remark}

\begin{proposition}\label{prop-polystable}\cite[Proposition 3.2]{K} Let $\sigma \in
\Sigma(Q,\beta)$ and $W \in \Rep(Q,\beta)^{ss}_{\sigma}.$ Then the
following are equivalent:
\begin{enumerate}
\renewcommand{\theenumi}{\arabic{enumi}}

\item $W$ is a direct sum of $\sigma$-stable representations;

\item $W$ is $\sigma$-polystable.

\end{enumerate}
\end{proposition}

\begin{remark} \label{rmk-polystable} Let $\sigma \in \Sigma(Q,\beta)$ and
consider the action of $\GL(\beta)$ on the variety
$\Rep(Q,\beta)^{ss}_{\sigma}$. Then for every $W \in
\Rep(Q,\beta)^{ss}_{\sigma}$ there exists a unique (up to
isomorphism) $\widetilde{W} \in \overline{\GL(\beta)W} \bigcap
\Rep(Q, \beta)^{ss}_{\sigma}$ such that $\widetilde{W}$ is
$\sigma$-polystable. More precisely, consider a Jordan-Holder
filtration of $W$ (in the category $\Rep(Q)^{ss}_{\sigma}$):
$$
F_{\bullet}(W): 0=W_0 \varsubsetneq W_1 \varsubsetneq \dots
\varsubsetneq W_{l-1} \varsubsetneq W_l=W,
$$
with $W_i / W_{i-1}$ $\sigma$-stable representations. Viewing the
associated graded representation $gr(F_{\bullet}(W))$ as the limit
of a $1$-parameter subgroup of $\GL(\beta)$ acting on $W$, we
deduce that
$$
gr(F_{\bullet}(W)):=\bigoplus_{i=1}^l W_i / W_{i-1} \in
\overline{\GL(\beta)W} \bigcap \Rep(Q, \beta)^{ss}_{\sigma}.
$$
(Technically speaking, $gr(F_{\bullet}(W))$ is isomorphic to a
point in the closure of the orbit of $W$.) In any case,
Proposition \ref{prop-polystable} implies that
$$
gr(F_{\bullet}(W))\cong \widetilde{W}.
$$
\end{remark}



The following lemma will come in handy when finding the faces of
GIT-cones (Proposition \ref{prop-essen1}):

\begin{lemma}\label{handy-lemma} Let $\sigma_1, \sigma_2 \in \Sigma(Q,\beta).$
\begin{enumerate}
\renewcommand{\theenumi}{\arabic{enumi}}

\item Assume that every $\sigma_1$-polystable representation
belongs to $\Rep(Q,\beta)^{ss}_{\sigma_2}.$ Then:
$$
\Rep(Q,\beta)^{ss}_{\sigma_1} \subseteq
\Rep(Q,\beta)^{ss}_{\sigma_2}.
$$

\item \cite[Lemma 4.1]{Re} Assume that
$$
\Rep(Q,\beta)^{ss}_{\sigma_1} \subseteq
\Rep(Q,\beta)^{ss}_{\sigma_2}.
$$

Let $W_2$ be $\sigma_2$-polystable. Then there exists a
$\sigma_1$-polystable representation $W_1$ such that $W_2 \in
\overline{\GL(\beta)W_1}.$ Moreover, $\Rep(Q,\beta)^{s}_{\sigma_2}
\subseteq \Rep^s_{\sigma_1}.$
\end{enumerate}
\end{lemma}

\begin{proof} $(1)$ Let $W \in \Rep(Q,\beta)^{ss}_{\sigma_1}$ and
let $\widetilde{W}$ be a $\sigma_1$-polystable representation so
that
$$\widetilde{W} \in \overline{\GL(\beta)W} \bigcap \Rep(Q,
\beta)^{ss}_{\sigma_1}.$$ Then, by assumption, $\widetilde{W}$ is
$\sigma_2$-semi-stable. But this clearly implies that $W$ is
$\sigma_2$-semi-stable, too. Indeed, one can use the description
\cite[Proposition 3.1]{K} of semi-stable representations in terms
of semi-invariants. So, $\Rep(Q,\beta)^{ss}_{\sigma_1} \subseteq
\Rep(Q,\beta)^{ss}_{\sigma_2}.$

$(2)$ See \cite[Lemma 4.1]{Re}.
\end{proof}

\section{Orbit cones of quiver representations} Fix a quiver $Q$ without oriented
cycles and a dimension vector $\beta \in \ZZ_{\geq 0}^{Q_0}.$

\begin{definition} Let $W \in \Rep(Q, \beta).$ The \emph{orbit cone} of $W$ is defined by
$$
\Omega(W)=\{\sigma \in \RR^{Q_0} \mid W \in
\Rep(Q,\beta)^{ss}_{\sigma}\}.
$$
\end{definition}

\begin{remark} Note that $\Omega(W)$ is a rational convex
polyhedral cone:
$$
\Omega(W)=\{\sigma \in \RR^{Q_0} \mid \sigma(\beta)=0 \text{~and~}
\sigma(\underline{d}_{W'})\leq 0, \text{~for all
subrepresentations~} W' \subseteq W \}.
$$
Furthermore, since there are finitely many dimension vectors
$\beta'$ with $\beta' \leq \beta,$ it follows that there are
finitely many orbit cones.
\end{remark}

First, a simple lemma:

\begin{lemma}\label{easy-lemma} Let $\sigma \in C(Q,\beta)$ and $W \in \Rep(Q,\beta)^{ss}_{\sigma}$ be
a $\sigma$-polystable representation. Then
$$
\sigma \in \rel(\Omega(W)).
$$
\end{lemma}

\begin{proof} Let $W'\subseteq W$ be a subrepresentation such
that $\sigma(\underline{d}_{W'})=0.$ This implies that $W',$
$W/W',$ and hence $W'\bigoplus W/W'$ are $\sigma$-semi-stable.
Now, it is well-known that $W'\bigoplus W/W' \in
\overline{\GL(\beta)W}.$ As the orbit of $W$ is closed in
$\Rep(Q,\beta)^{ss}_{\sigma},$ it follows
$$
W'\bigoplus W/W' \cong W.
$$

From this we obtain that $\Omega(W)=\Omega(W) \bigcap \mathbb
H(\underline{d}_{W'}).$ In other words, $\sigma$ can not lie on
any proper face of $\Omega(W)$ and this finishes the proof.
\end{proof}

Now, let us identify some of the faces of $\Omega(W)$ (compare to
\cite[Lemma 3.5]{Re}):

\begin{lemma}\label{lemma-face-stab} Let $\sigma_0 \in \Sigma(Q,\beta)$ and $W\in \Rep(Q,
\beta)^{ss}_{\sigma_0}.$ Consider a Jordan-H$\ddot{o}$lder
filtration of $W$ (in the category $\Rep(Q)^{ss}_{\sigma_0}$):
$$
F_{\bullet}(W): 0=W_0 \varsubsetneq W_1 \varsubsetneq \dots
\varsubsetneq W_{l-1} \varsubsetneq W_l=W,
$$
with $W_i / W_{i-1}$ $\sigma_0$-stable representations. Let
$\widetilde{W} \in \overline{\GL(\beta)W} \bigcap \Rep(Q,
\beta)^{ss}_{\sigma_0}$ be so that $\widetilde{W}$ is
$\sigma_0$-polystable. Then
$$
\Omega(\widetilde{W})=\Omega(W) \bigcap \{\sigma \in \RR^{Q_0}
\mid \sigma(\beta_i)=0, 1 \leq i \leq l \},
$$
where $\beta_i=\underline{d}_{W_i}, 1 \leq i \leq l.$ In
particular, $\Omega(\widetilde{W})$ is a face of $\Omega(W).$
\end{lemma}

\begin{proof} We have seen in Remark \ref{rmk-polystable} that
\begin{eqnarray}\label{eqn-gr}
gr(F_{\bullet}(W)):=\bigoplus_{i=1}^l W_i / W_{i-1} \cong
\widetilde{W}.
\end{eqnarray}

Now, if $\sigma \in \Omega(W) \bigcap \{\sigma \in \RR^{Q_0} \mid
\sigma(\beta_i)=0, 1 \leq i \leq l \}$ then $W_i$ are clearly
$\sigma$-semi-stable. This implies that $W_i / W_{i-1}$ are
$\sigma$-semi-stable and hence $\widetilde{W}$ is
$\sigma$-semi-stable.

Conversely, let $\sigma \in \Omega(\widetilde{W}).$ Then
$(\ref{eqn-gr})$ implies that the $W_i/W_{i-1}$ are
$\sigma$-semi-stable. Hence, $W=W_l$ is $\sigma$-semi-stable and
$\sigma(\beta_i)=0$ for all $1 \leq i \leq l.$ This finishes the
proof.
\end{proof}

The next proposition is essential for the description of the
GIT-classes and for the construction of the GIT-fan associated to
$(Q,\beta)$. It was proved for points of normal projective
varieties in \cite[Proposition 3.6]{Re}. The proof in that case
makes use of a fundamental finiteness result of Dolgachev and Hu
\cite[Theorem 1.5]{DH}. Our arguments are conceptually much
simpler and do not use Dolgachev and Hu's finiteness result.

\begin{proposition}[Faces of orbit cones]\label{prop-essen} Let $W \in \Rep(Q, \beta).$
Then the following are true.
\begin{enumerate}
\renewcommand{\theenumi}{\arabic{enumi}}

\item The faces of $\Omega(W)$ are exactly those of the form
$\Omega(W')$ with $W' \in \overline{\GL(\beta)W}.$

\item There exists $W_0 \in \overline{\GL(\beta)W}$ such that
$\Omega(W)=\Omega(W_0),$ and
$$
\sigma \in \rel(\Omega(W)) \Longleftrightarrow W_0 \text{~is~}
\sigma-\text{polystable}.
$$
\end{enumerate}
\end{proposition}

\begin{proof}$(1)$ First, let us show that if $W' \in
\overline{\GL(\beta)W}$ then $\Omega(W')$ is a face of
$\Omega(W)$. Assume $\Omega(W')$ is not the trivial cone (i.e., it
has positive dimension) and pick $\sigma' \in \rel(\Omega(W'))
\bigcap \ZZ^{Q_0}$. Let $W'' \in \overline{\GL(\beta)W'}\bigcap
\Rep(Q,\beta)^{ss}_{\sigma'}$ be so that $W''$ is
$\sigma'$-polystable. By Lemma \ref{lemma-face-stab}, we know that
$\Omega(W'')$ is a face of $\Omega(W')$. But $\Omega(W'')$
contains a relative interior point of $\Omega(W')$ and therefore,
we must have $\Omega(W'')=\Omega(W')$.

We clearly have that $W'' \in \overline{\GL(\beta)W} \bigcap
\Rep(Q,\beta)^{ss}_{\sigma'}$. Applying Lemma
\ref{lemma-face-stab} again, we obtain that $\Omega(W'')$ is a
face of $\Omega(W)$ and so $\Omega(W')$ is indeed a face of
$\Omega(W)$.

Now, let $\mathcal F$ be a proper face of $\Omega(W)$. If
$\mathcal F$ is the trivial face then $\mathcal F=\Omega(S)$,
where $S$ is the zero element of the vector space $\Rep(Q,\beta)$
(that is to say, $S$ is the unique semi-simple representation
belonging to $\Rep(Q,\beta)$). Now, let us assume that $\mathcal
F$ is not the trivial face and choose $\widetilde{\sigma} \in
\rel{\mathcal F}\bigcap \ZZ^{Q_0}$. Then, $\mathcal F$ is the
intersection of all facets of $\Omega(W)$ containing
$\widetilde{\sigma}$ and so we can write:
$$
\mathcal F=\Omega(W)\bigcap \bigcap_{W'} \{\sigma \in \RR^{Q_0}
\mid \sigma(\underline{d}_{W'})=0\},
$$
where the intersection on the right is over all subrepresentations
$W'\subseteq W$ with $\widetilde{\sigma}(\underline{d}_{W'})=0$.

Consider a Jordan-Holder filtration of $W$ (in the category
$\Rep(Q)^{ss}_{\widetilde{\sigma}}$):
$$
F_{\bullet}(W): 0=W_0 \varsubsetneq W_1 \varsubsetneq \dots
\varsubsetneq W_{l-1} \varsubsetneq W_l=W,
$$
with $W_i / W_{i-1}$ $\widetilde{\sigma}$-stable representations.
Let $\widetilde{W} \in \overline{\GL(\beta)W}\bigcap \Rep(Q,
\beta)^{ss}_{\widetilde{\sigma}}$ be so that $\widetilde{W}$ is
$\widetilde{\sigma}$-polystable. In what follows, we show that
$$
\mathcal F=\Omega(\widetilde{W}).
$$
We have seen in Lemma \ref{lemma-face-stab} that
$$
\Omega(\widetilde{W})=\Omega(W)\bigcap \bigcap_{i=1}^l \mathbb
H(\underline{d}_{W_i})
$$
and $\Omega(\widetilde{W})$ is a face of $\Omega(W)$. So,
$\mathcal F \subseteq \Omega(\widetilde{W})$ as
$\Omega(\widetilde{W})$ is a face containing a relative interior
point of $\mathcal F$.

Now, if $W'\subseteq W$ is a subrepresentation such that
$\widetilde{\sigma}(\underline{d}_{W'})=0$ then $W'$ becomes
$\widetilde{\sigma}$-semi-stable. The uniqueness of the factors of
the Jordan-Holder filtration $F_{\bullet}(W)$ implies that
$\underline{d}_{W'}$ is a non-negative integer linear combination
of vectors of the form
$\underline{d}_{W_i}-\underline{d}_{W_{i-1}}$. So,
$\sigma(\underline{d}_{W'})=0$ for all $\sigma \in
\Omega(\widetilde{W})$ and hence $\Omega(\widetilde{W})\subseteq
\mathcal{F}$. We have just proved that $\mathcal
F=\Omega(\widetilde{W})$ with $\widetilde{W} \in
\overline{\GL(\beta)W}$.

$(2)$ Let $\sigma_0 \in \rel(\Omega(W)) \bigcap \ZZ^{Q_0}$ and let
$W_0 \in \overline{\GL(\beta)W} \bigcap \Rep(Q,
\beta)^{ss}_{\sigma_{0}}$ be so that $W_0$ is
$\sigma_0$-polystable. We are going to show that $W_0$ satisfies
the required properties. Using $(1)$, it is clear that
$$
\Omega(W_0)=\Omega(W).
$$

The implication $"\Longleftarrow"$ follows immediately from Lemma
\ref{easy-lemma}. For the implication $"\Longrightarrow",$ let
$\sigma \in \rel(\Omega(W))$ and assume that $\GL(\beta)W_0$ is
not closed in $\Rep(Q, \beta)^{ss}_{\sigma}$. Then there exists
$W'' \in (\overline{\GL(\beta)W_0}-\GL(\beta)W_0)\bigcap \Rep(Q,
\beta)^{ss}_{\sigma}$. Note that $\sigma_0 \notin \Omega(W'')$ as
$\GL(\beta)W_0$ is closed in $\Rep(Q,\beta)^{ss}_{\sigma_{0}}.$
From this and part $(1)$, we get that  $\Omega(W'')$ is a proper
face of $\Omega(W_0)$. But this contradicts the fact that
$\Omega(W'')$ contains the relative interior point $\sigma$ of
$\Omega(W).$
\end{proof}




\section{The GIT-fan}
First a simple lemma:
\begin{lemma}\label{GIT-cone-lemma} If $\sigma \in C(Q,\beta)$ then
$$
C(\sigma)=\bigcap_{W} \Omega(W),
$$
where the intersection is over all $\sigma$-polystable
representations $W \in \Rep(Q,\beta).$
\end{lemma}

\begin{proof} By Corollary \ref{coro-inte-wt}, we can assume that
$\sigma \in \Sigma(Q,\beta).$ The inclusion $C(\sigma) \subseteq
\bigcap_{W} \Omega(W)$ is obvious while the other inclusion
follows immediately from Lemma \ref{handy-lemma}{(1)}. This
finishes the proof.
\end{proof}

It is clear that in order to prove Theorem \ref{thm-GIT-fan}, we
need to find the faces and relative interiors of GIT-cones.

\begin{proposition}[Faces of GIT-cones] \label{prop-essen1} Let $\sigma \in
C(Q,\beta) \bigcap \ZZ^{Q_0}$ and let $\mathcal F$ be a face of
$C(\sigma)$. If $\sigma_0 \in \rel(\mathcal F)\bigcap \ZZ^{Q_0}$
then
$$
\mathcal F = C(\sigma_0).
$$
In particular, any $\sigma_0 \in \rel(C(\sigma))\bigcap \ZZ^{Q_0}$
is GIT-equivalent to $\sigma$.
\end{proposition}

\begin{proof} By Lemma \ref{GIT-cone-lemma}, we know that
$$
C(\sigma)=\bigcap_{W} \Omega(W),
$$
where the intersection is over all $\sigma$-polystable
representations $W \in \Rep(Q,\beta)$. Moreover, this is a finite
intersection and so we have
$$
\mathcal F=\bigcap_{W} \mathcal F_{W,\sigma_0},
$$
where $\mathcal F_{W,\sigma_0}$ is the unique face of $\Omega(W)$
such that $\sigma_0 \in \rel(\mathcal F_{W,\sigma_0})$.

Now, for every $\sigma$-polystable representation $W \in
\Rep(Q,\beta),$ let $\widetilde{W} \in \Rep(Q,\beta)$ be the
unique (up to isomorphism) $\sigma_0$-polystable representation
such that $\widetilde{W} \in \overline{\GL(\beta)W}$. From the
proof of Proposition \ref{prop-essen}, we know that each $\mathcal
F_{W,\sigma_0}=\Omega(\widetilde{W})$ and hence
$$
\mathcal F=\bigcap_{W}\Omega(\widetilde{W}),
$$
where the intersection is over all $\sigma$-polystable
representations $W \in \Rep(Q,\beta)$. Using Lemma
\ref{handy-lemma}{(2)} and Lemma \ref{GIT-cone-lemma}, this
intersection is exactly $C(\sigma_0)$ and so  $\mathcal
F=C(\sigma_0)$.
\end{proof}





The next result gives a polyhedral description of a GIT-class
(compare to \cite[Lemma 4.2 \footnote{Note that Rassayre's proof
works only for rational points of the ample cone.}]{Re}):

\begin{proposition}[GIT-classes]\label{prop-GIT-classes} Let $\sigma \in C(Q,\beta)$ and $F=\langle
\sigma \rangle$. Then
$$
F=\rel(C(\sigma)).
$$
\end{proposition}

\begin{proof} First of all, we can assume that $\sigma \in
\Sigma(Q,\beta)$ by Corollary \ref{coro-inte-wt}. Now, let
$\sigma' \in \rel(C(\sigma))$.

As $\rel(C(\sigma'))\bigcap \rel(C(\sigma)) \neq \emptyset$, we
know that $C(\sigma')\bigcap C(\sigma)$ is a rational convex
polyhedral cone whose relative interior is non-empty. Hence, we
can choose $$\sigma_0 \in \rel(C(\sigma'))\bigcap \rel(C(\sigma))
\bigcap \ZZ^{Q_0}.$$ This implies that $\sigma_0 \thicksim \sigma$
by Proposition \ref{prop-essen1}. We also have that
$\Rep(Q,\beta)^{ss}_{\sigma} \subseteq
\Rep(Q,\beta)^{ss}_{\sigma'} \subseteq
\Rep(Q,\beta)^{ss}_{\sigma_0}$ and so $\sigma' \thicksim \sigma$.
We have jus proved
$$
\rel(C(\sigma)) \subseteq F.
$$
Now, let us prove the other inclusion. From Lemma
\ref{easy-lemma}, we know that
$$
\sigma \in \bigcap_{W} \rel(\Omega(W)),
$$
where the intersection is over all $\sigma$-polystable
representations $W \in \Rep(Q,\beta).$ Since this is a non-empty
intersection and we work with a finite intersection of rational
convex polyhedral cones, we have:
$$
\bigcap_{W} \rel(\Omega(W))=\rel(\bigcap_W \Omega(W)).
$$
Using Lemma \ref{GIT-cone-lemma} and applying Lemma
\ref{easy-lemma} for every point of $F,$ we obtain $F \subseteq
\rel(C(\sigma)).$



\end{proof}

\begin{proof}[Proof of Theorem \ref{thm-GIT-fan}]
$(1)$ This part follows now from Lemma \ref{finiteness-lemma} and
Proposition \ref{prop-GIT-classes}.

$(2)$ It is clear that $\mathcal F(Q,\beta)=\{C(\sigma) \mid
\sigma \in C(Q,\beta) \}$ is a finite cover of $C(Q,\beta)$.
Moreover, Proposition \ref{prop-essen1} shows that every face of a
GIT-cone is a again a GIT-cone. It remains to prove that the
intersection of two GIT-cones is also a GIT-cone. Let $\sigma_1,
\sigma_2 \in C(Q,\beta)$ be so that $C(\sigma_1) \bigcap
C(\sigma_2) \varsupsetneq \{0\}$. Now, pick a $\sigma_0 \in
\rel(C(\sigma_1)\bigcap C(\sigma_2)) \bigcap \ZZ^{Q_0}$. Let $F_i$
be the unique (not necessarily proper) face of $C(\sigma_i)$ so
that $\sigma_0 \in \rel(F_i),$ $i \in \{1,2\}$. From Proposition
\ref{prop-essen1} follows $F_1=F_2=C(\sigma_0)$. So, $C(\sigma_0)$
is a face of $C(\sigma_1)\bigcap C(\sigma_2)$ containing a
relative interior point this intersection. Therefore,
$$
C(\sigma_1) \bigcap C(\sigma_2)=C(\sigma_0)
$$
and this finishes the proof.
\end{proof}

\end{document}